\documentclass[openacc]{rsproca_new}

\usepackage{silence}
\usepackage[noabbrev,nameinlink]{cleveref}

\usepackage[caption=false]{subfig}

\definecolor{red}{RGB}{207,18,61}
\definecolor{blue}{RGB}{0,137,189}

\newcommand\ie{i.e.\ }

\newcommand{\ts}{\textsuperscript}

\newcommand{\etal}{\textit{et al.}}

\newcommand{\Rey}{\mbox{\textit{Re}}}
\newcommand{\Ca}{\mbox{\textit{Ca}}}

\newcommand*{\tran}{^{\mkern-1.5mu\mathsf{T}}} 
\renewcommand{\matrix}[1]{
  \begin{pmatrix}
    #1
  \end{pmatrix}
} 

\makeatletter
\newcommand*{\diff}%
    {\@ifnextchar^{\DIfF}{\DIfF^{}}}
\def\DIfF^#1{%
    \mathop{\mathrm{\mathstrut d}}%
        \nolimits^{#1}\gobblespace}
\def\gobblespace{\futurelet\diffarg\opspace}
\def\opspace{%
    \let\DiffSpace\!%
    \ifx\diffarg(%
        \let\DiffSpace\relax
    \else
        \ifx\diffarg[%
            \let\DiffSpace\relax
        \else
            \ifx\diffarg\{%
                \let\DiffSpace\relax
            \fi\fi\fi\DiffSpace}

\titlehead{Research}
\jname{rspa}
\Journal{Proc R Soc A\ }

\title{Non-linear estimators for the observation and stabilisation of falling liquid films}
\author{
  Oscar Holroyd\(^{1}\),
  Radu Cimpeanu\(^{1}\), and
  Susana N. Gomes\(^{1}\)
}
\address{\(^{1}\)
  Warwick Mathematics Institute,
  University of Warwick,
  Coventry,
  CV4 7AL,
  UK
}
\corres{Susana N. Gomes\\\email{susana.gomes@warwick.ac.uk}}
\subject{applied mathematics, control theory, fluid dynamics}

\begin{document}
  \begin{abstract}
    Falling liquid films are complex physical processes modelled by infinite-dimensional dynamical systems. The interfacial behaviour is characterised by the Reynolds number, a nondimensional parameter, and above a critical value the uniform film becomes unstable and travelling waves occur. By injecting and removing fluid from the base at discrete locations, we aim to stabilise otherwise unstable flat interfaces, with the additional limitation that observations of the state are restricted to finitely many measurements of the film height. Successful feedback control has recently been achieved in the case of full observations using linear-quadratic regulator controls coupled to asymptotic approximations of the Navier-Stokes equations, but restricted observations severely curtailed their performance. In this study we couple the well understood full-information feedback control strategy to a non-linear estimator. The dynamics of the estimator are designed to approximate those of the film, and we apply a forcing term chosen to ensure that measurements of the estimator match the available measurements of the film. Using this method, we restore the performance of the controls to a level approaching their full information counterparts, even at moderately large Reynolds numbers. We also briefly investigate the effects of the relative positioning of actuators and observers on the resulting dynamics.
  \end{abstract}
  \maketitle

  \section{Introduction}%
  \label{sec:introduction}
  The observation and stabilisation of thin liquid films is a prototypical problem at the intersection of fluid dynamics, control theory, and asymptotic analysis. The system itself displays a rich and varied range of parameter-dependent features: flat films, travelling waves, and chaotic behaviour. This high degree of complexity emerging from a relatively simple system makes it a prime candidate for exploring the developing field of infinite-dimensional non-linear control. There are numerous other systems with control applications exhibiting complex behaviour -- multiphase flows in ice buildup on aerofoils~\cite{palacios2011instantaneous}, obstacle avoidance in self-driving vehicles~\cite{paden2016survey}, and crowd management~\cite{burger2013mean} to name but a few. Additionally, in most of these situations it is either impractical or impossible to perfectly observe the dynamics.

  Falling liquid films provide a balance between model complexity and analytical tractability of the underlying fluid-dynamics to allow a focus on bridging the gap between non-linear partial differential equations (PDEs) and control theory. The flow of a thin liquid film, falling under gravity, becomes unstable above a critical Reynolds number. This begins with the development of two-dimensional (2D) waves, and progresses to three-dimensional (3D) spatiotemporal chaos. This was demonstrated experimentally by Nosoko \etal~\cite{nosoko1996characteristics}, and numerically by Papageorgiou and Smyrlis~\cite{papageorgiou1991route}. Most relevant to this work, a number of thin-film asymptotic models -- where the height of the film is assumed to be much smaller than the characteristic wavelength -- have been developed. Many of these are covered extensively in a pair of reviews by Kalliadasis \etal~\cite{kalliadasis2011falling} and Craster and Matar~\cite{craster2009dynamics}. Such models range from relatively simple equations for the dynamics of the film height to more complex systems connecting film height, downstream flux, and other more abstract quantities. Examples of the first type include the long-wave expansion by Benney~\cite{benney1966long} from the 1960s, perhaps the original falling liquid film PDE, and Ooshida's mixed-derivative regularised equation~\cite{ooshida1999surface}. Ruyer-Quil and Manneville~\cite{ruyer2000improved,ruyer2002further} pioneered a new style of weighted-residual models in the 2000s, which has led to more recent work by Richard \etal~\cite{richard2019optimization}, Usha \etal~\cite{usha2020evolution}, and Mukhopadhyay \etal~\cite{mukhopadhyay2022modelling}. These studies have provided insight into how two- and four-equation models result in far superior accuracy that for previous techniques, raising the upper limit of the Reynolds numbers to which the thin-film assumption remains relevant. Ruyer-Quil \etal~\cite{ruyer2014dynamics} gives an overview of some of these low-dimensional modelling attempts (notably missing some of the more recent work).

  Applications of controls to falling liquid films include coating flows~\cite{weinstein2004coating}, heat and mass transfer~\cite{wayner1991effect} and thermoelectric micro-channel cooling~\cite{darabi2001electrohydrodynamic}. Some cases, such as coating flows, require flat interfaces, whereas for other uses -- like heat transfer -- the larger surface area provided by highly corrugated films is more desirable. Early investigation of control techniques focussed principally on the physics of potential actuation mechanisms, largely limited to passive or simple controls. Unsurprisingly, the large number of applications has resulted in a wide range of possible external inputs, and so there is an extensive body of literature studying how these different inputs impact the stability and the critical Reynolds number. These include heating and cooling of the fluid~\cite{bankoff1971stability}, electric~\cite{tseluiko2006wave} or magnetic~\cite{amaouche2013hydromagnetic} fields, porous~\cite{ogden2011gravity}, deformable~\cite{gaurav2007stability}, or corrugated~\cite{dull2025spatio} substrates, and many more. In this work we focus on blowing and suction through the base by a small number of discrete actuators~\cite{thompson2016falling}. Despite difficulties in its real-world instantiation, such controls (restricted to discrete portions of the boundary) are an improvement over either distributed controls commonly used in theoretical studies of non-linear PDE control~\cite{chow2016control} or continuous boundary controls~\cite{aljamal2018linearized}, which, for the blowing and suction controls considered here, would be impossible to construct experimentally, even though they are more convenient analytically and computationally. Additionally, they act on the film at the same scale as the perturbations (unlike air-jets, which can be too strong~\cite{lhuissier2016blowing,ojiako2020deformation}, or surfactants, which can be too weak to use as a control~\cite{blyth2004effect}), which makes them more computationally tractable due to the reduced time of action which results in shorter simulations and more efficient control.

  More recent literature has extended beyond investigating control mechanisms to how these methods can be used to implement closed-loop control schemes, where information from the system is instantaneously fed-back in as part of a control term. To start with, work focused on the Kuramoto-Sivashinsky (KS) equation, a highly simplified thin-film model that is the simplest PDE to exhibit chaotic behaviour~\cite{brummitt2009search}. Armaou and Christofides~\cite{armaou2000feedback,christofides2000global} and more recently Gomes \etal~\cite{gomes2015controlling,gomes2017controlling} applied control theoretical techniques (see~\cite{zabczyk2020mathematical} for an introduction) from both analytical and computational perspectives to show the potential for successful results in this space. Thompson \etal~\cite{thompson2016stabilising} extended these results to Benney and weighted-residual systems, although the increased complexity of these models precluded any proof of convergence to the steady state. Cimpeanu \etal~\cite{cimpeanu2021active} were the first to attempt, with some success, to apply controls to the full thin-film system, modelled by the Navier-Stokes (albeit using the simple proportional controls used by Thompson \etal~\cite{thompson2016falling}). In the current authors' previous work they combined the two approaches to show that a chain of successively simplified approximations can be used to design linear-quadratic regulator (LQR) controls that are able to control simulations of the Navier-Stokes film~\cite{holroyd2024linear}.

  Although there is a rich body of literature on the control of infinite dimensional systems~\cite{curtain2012introduction}, much of it deals only with either linear problems~\cite{zwart2020optimal,bergeling2020closed}, or simple non-linear problems with full-information~\cite{capistrano2023rapid}. In particular, most previous work on thin-film control are focussed solely on control and assume that perfect information of the system is available. This is, of course, trivial when the control system is being simulated on a computer, but in reality can be very challenging~\cite{mouza2000measurement,ishikawa2004measurement}. In a previous study~\cite{holroyd2024stabilisation}, we repurposed the Luenberger observer/estimator~\cite{luenberger1971introduction} approach used for the reduced-dimension problem~\cite{thompson2016stabilising} for use with direct numerical simulations (DNS) of the Navier-Stokes equation. Despite a narrowing of the region of the parameter space in which the controls were successful, this did show that thin-film control is possible with severely limited information.

  After building the necessary background on the fundamental thin-film problem in \cref{sec:the_thin_film_problem}, and the basics of the control problem in \cref{sec:stabilising_the_flat_film}, we move on to a discussion of the techniques available to tackle the issue of restricted observations in \cref{sec:restricted_observations}, with recent contributions in the past years from our own work and in the wider research area acting as foundational elements. There we cover the two methods used in previous work before expanding into the new, more complex estimator design in \cref{sub:observation_and_estimation}. Finally, in \cref{sub:estimation_and_control} we show how this can be coupled with the standard LQR controller to stabilise the flat interfacial profile in more challenging areas of the parameter space. We also briefly explore the interplay between actuator and observer locations, which provides the problem with an additional degree of complexity.

  \section{The Thin Film Problem}%
\begin{figure}[!htb]
    \centering
    \includegraphics[width=\textwidth]{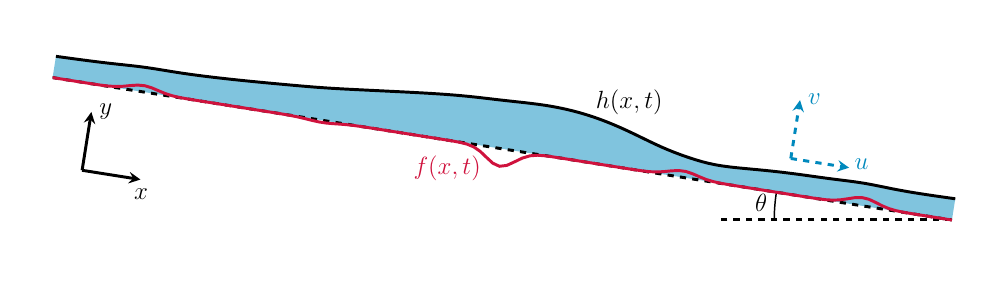}
    \caption{Diagram of the falling liquid film control problem. A thin liquid film with height \(h(x,t)\) falls under gravity down a plane inclined at an angle \(\theta\). Stabilisation of the interface occurs via a basal forcing \(f\) applied at the lower boundary \(y = 0\). The forcing consists of a finite number of actuators, injecting or removing fluid from the system.
    }%
    \label{fig:thin_film}
  \end{figure}

  \label{sec:the_thin_film_problem}
  As shown in \cref{fig:thin_film}, we consider the flow of a thin film of liquid falling down an inclined plane, tilted at an angle \(0 < \theta < \pi / 2\) from the horizontal. A combination of sufficiently weak surface tension or viscosity, strong inertial effects, and large inclination angle results in a gravity-induced destabilisation of the flat, uniform interface. Left unchecked, this results in the development of a non-linear wave, which travels down the plane. To combat growing, unstable perturbations to the film, we are allowed to inject or remove fluid from the otherwise impermeable, no-slip base plate at a finite number of evenly-spaced points, known as actuators.

  Mathematically we describe the problem in a periodic reference frame rotated to match the inclination angle, where the \(x\)- and \(y\)-directions correspond to the streamwise and normal vectors respectively. The interface between the liquid and the gas above it is described by a function \(h(x,t) > 0\). This implies that the travelling waves cannot break (which would require \(h\) to be multivalued) and that the film cannot rupture (which would mean that \(h(x,t) = 0\) for some value of \(x\) and \(t\)).

  The dynamics of the film is governed by the external gravitational forcing \(\vec{g} = (0, g)\), the parameters describing the domain -- the angle \(\theta\) and the dimensional length \(L^*\), and the physical properties of the fluid -- density \(\rho\), viscosity \(\mu\), and surface tension \(\gamma\). For liquid-gas systems such as those considered here, the density (and viscosity) of the gas tends to be several orders of magnitude lower than that of the liquid. This means that the inertial forces associated with the liquid film completely dominate that of the gas, and so we can ignore the gas dynamics entirely.

  Although we have thus far described the film in a dimensional setting, it is convenient to nondimensionalise the problem before we continue. For a film with total mass corresponding to a mean height \(h_N\), the uncontrolled system admits a uniform solution known as the Nusselt solution~\cite{nusselt1923warm}:
  \begin{equation}
    h(x, t) = h_N, \qquad u_N(x, t) = \frac{\rho g h_N^2 \sin\theta}{2\mu},
  \end{equation}
  where \(u_N\) is the horizontal velocity at the surface, \(y=h(x,t)\). Using the height of the film \(h_N\) as the length scale, the Nusselt velocity \(u_N\) as the velocity scale, and \(\frac{\rho u_N}{h_N}\) as the pressure scale (chosen to balance the pressure gradient \(p_x\) with viscous forces \(\mu u_{yy}\)~\cite{kalliadasis2011falling}), we can nondimensionalise the problem. This leaves us with the existing dimensionless parameter \(\theta\), the dimensionless domain length \(L = L^* / h_N\), and two new dimensionless numbers
  \begin{equation}
    \Rey = \frac{\rho u_N h_N}{\mu}, \qquad \Ca = \frac{\mu u_N}{\gamma}.
  \end{equation}
  The Reynolds number \(\Rey\) describes the relative importance of inertial and viscous effects, and the capillary number \(\Ca\) measures the relative importance of gravity and surface tension. Most thin-film physical regimes are captured by the dimensionless framework given by \(1 < \Rey < 200\) and \(0.001 < \Ca < 0.05\).

  \subsection{Navier-Stokes Equations}%
  \label{sub:navier_stokes}
  For the parameter regimes considered here where the fluid is viscous and incompressible, we can use the Navier-Stokes equations to describe the dynamics of the fluid flow. These equations are solved for the velocity \(\vec{u}(x, y, t) = (u(x, y, t), v(x, y, t))\) and pressure \(p(x, y, t)\). They consist of the momentum equations
  \begin{align}
    \label{eqn:momentum1}\Rey (u_t + uu_x + vu_y) &= -p_x + 2 + u_{xx} + u_{yy},\\
    \label{eqn:momentum2}\Rey (v_t + uv_x + vv_y) &= -p_y + 2\cot\theta + v_{xx} + v_{yy},
  \end{align}
  and the continuity equation
  \begin{equation}
    \label{eqn:continuity}u_x + v_y = 0.
  \end{equation}

  The PDE is closed with boundary conditions at the base and at the interface. At the base, \(y=0\), we have no-slip and no penetration other than that due to the actuators:
  \begin{equation}
    \label{eqn:base}u = 0, \qquad v = f.
  \end{equation}
  At the interface, \(y = h(x, t)\), we have the non-linear dynamic stress balance~\cite{kalliadasis2011falling}
  \begin{align}
    \label{eqn:stress1}(v_x + u_y)(1-h_x^2) + 2h_x (v_y  -u_x) &= 0, \\
    \label{eqn:stress2}p - \frac{2(v_y + u_x h_x^2 - h_x(v_x + u_y))}{1+h_x^2} &= -\frac{1}{\Ca}\frac{h_{xx}}{(1+h_x^2)^{3/2}}.
  \end{align}
  Finally we have the kinematic equation
  \begin{equation}
    \label{eqn:kinematic}h_t = v - uh_x,
  \end{equation}
  which describes how the interface is advected by the velocity field. Defining the down-slope flux
  \begin{equation}
    q(x,t) = \int_0^h u(x, y, t) \diff y,
  \end{equation}
  we can combine \cref{eqn:continuity,eqn:base,eqn:kinematic} to obtain the 1D integrated mass-conservation equation
  \begin{equation}
    \label{eqn:conservation}h_t + q_x = f.
  \end{equation}

  The Navier-Stokes equations are our most accurate representation of the fluid flow, and so, while insight from physical experiments would be invaluable, here we perform highly accurate computational experiments by solving the Navier-Stokes equations with the open-source volume-of-fluid code Basilisk~\cite{basilisk2025}, as a closely aligned theoretical counterpart of the real-world setup. More details of this framework can be found in \cref{sec:direct_numerical_simulation}. We highlight the promising recent development of similar types of feedback control being applicable to Hele-Shaw cells in experiments~\cite{ayoubi2025control}. The key practical differences in our present system being the non-confined nature of the target falling film flow, as well as the imposition of the active control inside the rapidly moving liquid.

\begin{figure}[!htb]
    \centering
    \includegraphics[width=\textwidth]{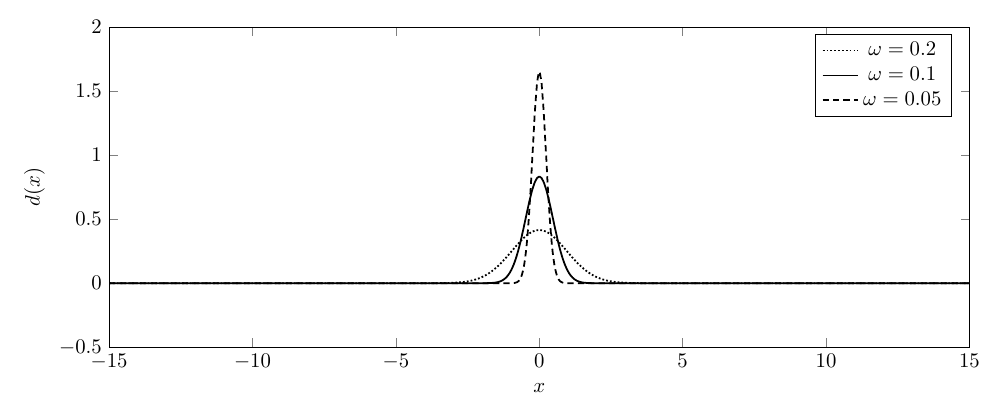}
    \caption{Actuator functions \(d(x)\) for varying values of the shape parameter \(\omega\). As \(\omega \to 0\), \(d\) converges to the Dirac delta \(\delta\).}%
    \label{fig:actuators}
  \end{figure}

  Since arbitrary injection and removal of fluid from the lower boundary is physically unfeasible, we restrict the forcing term to a finite number of injection sites (known as actuators):
  \begin{equation}
    f(x, t) = \sum_{i = 1}^{m} a_i(t) d(x - x_i),
  \end{equation}
  where \(a_i\) is the strength of the \(i\)\ts{th} actuator, which is located at \(x_i\). The individual actuation functions are
  \begin{equation}
    d(x) = \alpha \exp\left[ \frac{\cos(2\pi x/L) - 1}{\omega^2} \right].
  \end{equation}
  This is a smooth, periodic approximation to the Dirac delta distribution, chosen because although the Dirac delta is convenient from an analytical perspective (as used by Gomes \etal~in their proof of controllability for the KS equation~\cite{gomes2015controlling}), it is physically impossible and results in numerical instabilities in both the control design and in simulations. We also note that a smooth velocity profile is closer to what we might expect from the inlet in a physical experiment. The constant \(\alpha\) is chosen so that \(\int_0^L d(x) \diff x = 1\), and the width of the actuator is set by the parameter \(\omega\) -- we can see how the shape changes as \(\omega\) varies in \cref{fig:actuators}.

  \subsection{Benney Equation}%
  \label{sub:benney}
  Given the complexity of the Navier-Stokes equations \cref{eqn:momentum1,eqn:momentum2,eqn:continuity,eqn:base,eqn:stress1,eqn:stress2,eqn:kinematic}, it is unsurprising that the full, 2D system is challenging to control (we will discuss this in more detail in \cref{sec:stabilising_the_flat_film}). One particular reason is that the control term \(f\) acts only through the lower boundary, rather than being distributed throughout the domain~\cite{fattorini1968boundary}. Therefore, the first step towards control is to eliminate this problem.

  The first step in simplifying the problem is to introduce a long-wave parameter \(\epsilon = h_N / L^* = 1 / L\). To retain inertial and surface tension effects we assume that \(\Rey = O(1)\) and \(\Ca = O(\epsilon^2)\), and we match the magnitude of the control to that of the vertical velocity, \(f = O(\epsilon)\).

  One possible second step is to find an expression for the flux, \(q = Q(h, f)\), closing \cref{eqn:conservation}. By performing asymptotic expansions of \(u\) in powers of \(\epsilon\), Benney~\cite{benney1966long} derived a single-equation model for a falling liquid film. Extending this analysis to include the effects of \(O(\epsilon)\) controls~\cite{thompson2016falling}, we have
  \begin{equation}
    \label{eqn:benney}h_t = -\left[ \frac{h^3}{3}\left( 2 - 2 h_x \cot\theta + \frac{h_{xxx}}{\Ca} \right) + \Rey \left( \frac{8h^6h_x}{15} - \frac{2h^4f}{3} \right) \right]_x + f.
  \end{equation}

  Whilst the problem remains highly non-linear, the three, 2D, Navier-Stokes equations have been reduced down to a single, 1D, equation. More importantly, we have integrated out the upper and lower boundaries and are left with just the periodic boundary condition in the horizontal direction. This means that the control term \(f\) appears within the PDE, and its effects are distributed rather than acting locally at the boundary.

  Unfortunately, although the Benney equation is consistent at \(O(\epsilon)\)~\cite{ruyer2014dynamics}, and is asymptotically consistent as \(\Rey \to 0\), it only provides a reasonable approximation to the film evolution for small Reynolds numbers~\cite{scheid2005validity}. Above a flow-condition dependent critical Reynolds number a single-peak wave develops singular behaviour, which results in unphysical finite-time blow-up\cite{pumir1983solitary}. This can be seen in \cref{fig:comp_interface}.

  \subsection{Weighted-Residual Equations}%
  \label{sub:weighted_residual}
  To get around the unphysical behaviour inherent in single-equation thin film models, Ruyer-Quil and Manneville~\cite{ruyer2000improved,ruyer2002further} developed a weighted-residual methodology based on approximating the horizontal velocity as a sum of quadratic basis functions, designed to satisfy the no-slip boundary conditions at \(y=0\) and zero tangential stress at \(y=h\). Truncating at first order in \(\epsilon\), this allows for the integration of the system between the base and the interface, eliminating the vertical dimension while retaining separate dynamics for the flux \(q\).

  Extended to include the effects of forcing~\cite{thompson2016falling}, the momentum equations \cref{eqn:momentum1,eqn:momentum2} and normal and tangential stress conditions \cref{eqn:stress1,eqn:stress2} can be integrated to find an equation for the evolution of the flux:
  \begin{equation}
    \label{eqn:wr}\frac{2\Rey}{5}h^2q_t + q ={} \frac{h^3}{3}\left( 2 - 2 h_x \cot\theta + \frac{h_{xxx}}{\Ca} \right) + \Rey \left( \frac{18q^2h_x}{35} - \frac{34hqq_x}{35} + \frac{hqf}{5} \right).
  \end{equation}

\begin{figure}[!htb]
    \centering
    \includegraphics[width=\textwidth]{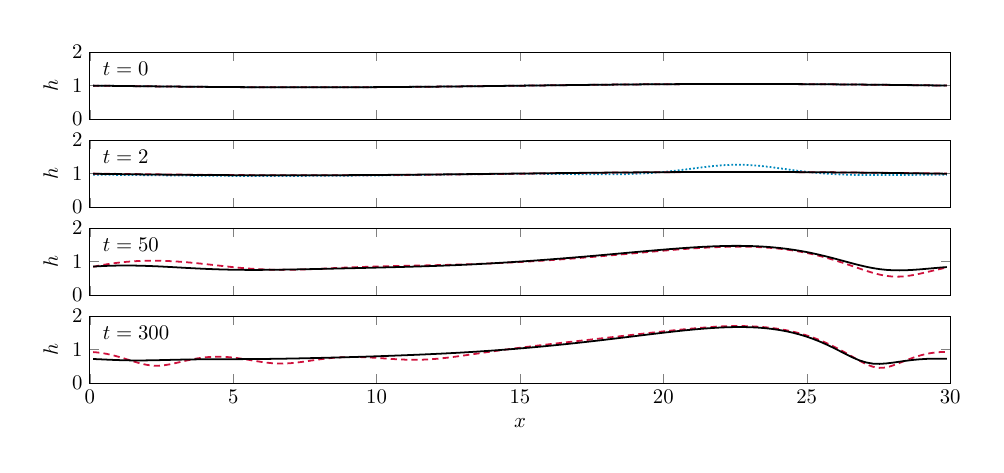}
    \caption{Comparison of the interfacial height \(h\) over time for Navier-Stokes, weighted residual, and Benney systems (black, red dashed, and blue dotted respectively) at \(\Rey = 10\), \(\Ca = 0.05\), \(\theta = \pi/3\). Although the weighted-residual model develops oversized capillary waves by \(t=300\), this is an improvement over the Benney equation, which blows up extremely rapidly after \(t = 2\) (too rapidly to easily capture). Peaks are shifted to \(22.5\) to allow easy comparison of interfacial shapes.}%
    \label{fig:comp_interface}
  \end{figure}

  Unlike the Benney equation, \cref{eqn:benney}, the weighted-residual system is well-behaved at all Reynolds numbers. As shown in \cref{fig:comp_interface}, it captures the behaviour of the full Navier-Stokes film with predictive accuracy at low-to-moderate values of \(\Rey\) (up to around \(\Rey = 7\)), but continues to exhibit qualitatively similar behaviour as \(\Rey\) continues to increase.

  \section{The Control Problem}%
  \label{sec:stabilising_the_flat_film}
  Thus far we have examined the problem from mathematical modelling and fluid dynamics perspectives, aiming to accurately approximate the dynamics of falling liquid films without considering the possibility of using the effects of the forcing from the base for control. We will now reframe the problem in a control setting.

  \subsection{General Control Problem}%
  \label{sub:general_control_problem}
  Consider the general control problem
  \begin{align}
    \label{eqn:cont_full_1}\dot{x} &= \mathcal{A}x + \mathcal{B}u, \\
    \label{eqn:cont_full_2}y &= \mathcal{C}x, \\
    x(0) &= x_0.
  \end{align}
  Here \(x\) is the state of the system (in an appropriately chosen functional Hilbert space \(\mathcal{X}\)) which we wish to control towards the target state \(x^*\). The differential operator \(\mathcal{A}: \mathcal{X} \to \mathcal{X}'\) (where \(\mathcal{X}'\) is another Hilbert space with potentially fewer derivatives) describes the dynamics of the uncontrolled system. The control input \(u \in \mathcal{U}\) influences the system via the operator \(\mathcal{B}: \mathcal{U} \to \mathcal{X}'\), and we are provided an observation of the state \(y \in \mathcal{Y}\) by the operator \(\mathcal{C}: \mathcal{X} \to \mathcal{Y}\). If the initial condition \(x_0\) is not within the basin of attraction of \(x^*\), we must construct the function \(u(t)\) so that \(x \to x^*\) as \(t \to \infty\).

  Solutions to an open control problem such as this one are strongly dependent on the initial condition, and thin films exhibit chaotic behaviour~\cite{kalliadasis2011falling}. We would therefore expect that open-loop controls of the form \(u(t; x_0)\) would be highly sensitive to noise, whether that be in the form of uncertainty in either the parameters or the initial condition in physical experiments, or floating-point and convergence errors in numerical experiments. To combat this we add the additional constraint that \(u\) must be a feedback control, \ie \(u = u(y(t))\), which means we can simplify the control problem to a single equation
  \begin{equation}
    \label{eqn:cont_compact}\dot{x} = \left(\mathcal{A} + \mathcal{B}\mathcal{K}\mathcal{C}\right)x,
  \end{equation}
  where \(\mathcal{K}: \mathcal{Y} \to \mathcal{U}\) is the gain operator which we need to find.

  For most solvable systems \((\mathcal{A}, \mathcal{B}, \mathcal{C})\) there will be an infinite number of solutions. Introducing the quadratic cost
  \begin{equation}
    \label{eqn:cont_cost}\kappa = \int_{0}^{\infty} \langle x, \mathcal{Q} x \rangle + \langle u, \mathcal{R} u \rangle \diff t,
  \end{equation}
  where \(\mathcal{Q}: \mathcal{X} \to \mathcal{X}\) and \(\mathcal{R}: \mathcal{U} \to \mathcal{U}\) are symmetric positive definite operators, we limit our search to the solution to \cref{eqn:cont_compact} that minimises \cref{eqn:cont_cost}.

  Unfortunately control methodologies for general, non-linear, infinite-dimensional problems remain undeveloped, with simple non-linear PDEs~\cite{cerpa2022carleman,capistrano2023rapid,afshar2023extended} and even linear PDEs~\cite{bergeling2020closed} being an area of active research. To make our problem tractable, we make some restrictive but necessary assumptions, namely that the spaces \(\mathcal{X}, \mathcal{Y}\), and \(\mathcal{U}\) are finite dimensional, and that the operators \(\mathcal{A}, \mathcal{B}, \mathcal{C}, \mathcal{Q}\), and \(\mathcal{R}\) are linear. We are thus left with the task of finding the matrix \(K\) that minimises the discrete cost
  \begin{equation}
    \label{eqn:cost}\kappa = \int_{0}^{\infty} \langle x, Q x \rangle + \langle u, R u \rangle \diff t,
  \end{equation}
  subject to the constraint
  \begin{equation}
    \label{eqn:compact}\dot{x} = \left(A + BKC\right)x.
  \end{equation}
  The stabilisation, observation, and control of systems of linear ordinary differential equations (ODEs) is a well-studied problem from both a theoretical and applied standpoint~\cite{zabczyk2020mathematical}, and robust methodologies have been developed for linear-quadratic regulator (LQR) problems (linear ODEs with a quadratic cost), albeit typically for systems with a small number of equations. Larger systems have been used to approximate the control of infinite-dimensional systems, but typically only for linear PDEs~\cite{zwart2020optimal}.

  \subsection{Control of Thin Films with Full Information}%
  \label{sub:control_with_full_information}
  It is clear that the approach outlined above is not applicable to the Navier-Stokes problem in and of itself. \Cref{eqn:cost,eqn:compact} are linear and finite dimensional, whereas the Navier-Stokes equations are non-linear and infinite-dimensional. Furthermore, \cref{eqn:cont_full_1} describes a distributed control problem, and the Navier-Stokes control acts through a boundary condition. We have seen how integrating out the vertical component to derive the Benney and weighted-residual equations can solve the latter problem without losing significant levels of information, especially in the weighted-residuals case.

\begin{figure}[!htb]
    \centering
    \includegraphics[width=\textwidth]{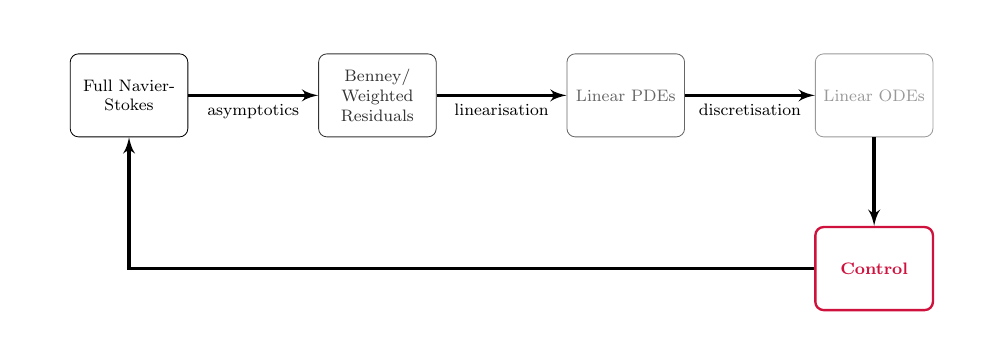}
    \caption{Diagrammatic illustration of the hierarchy of simplifications that permit the design of stabilising controls for DNS of a Navier-Stokes falling film. By making a chain of assumptions that increasingly abstract the model from the original system we eventually reach a point where established control theoretical results can be used. These controls can then be applied to observations of the full system projected down to the discrete space.}%
    \label{fig:hierarchy}
  \end{figure}

  To reach the form of \cref{eqn:cost,eqn:compact} we need to continue to simplify the thin film models \cref{eqn:benney} and \cref{eqn:conservation,eqn:wr}. We do this by first linearising the two PDEs (so \(h \approx 1 + \hat{h}\) and \(q \approx 2/3 + \hat{q}\)) to get
  \begin{equation}
    \label{eqn:lin_benney}\hat{h}_t = \left[ -2\partial_x + \left( \frac{2\cot{\theta}}{3} - \frac{8\Rey}{15} \right)\partial_{xx} - \frac{1}{3\Ca}\partial_{xxxx} \right]\hat{h} + \left[ 1 + \frac{2\Rey}{3}\partial_x \right] \hat{f},
  \end{equation}
  and
  \begin{align}
    \label{eqn:lin_wr1}\hat{h}_t &= -\hat{q}_x + \hat{f}, \\
    \label{eqn:lin_wr2}\hat{q}_t &= \left[ \frac{5}{\Rey} + \left( \frac{4}{7} - \frac{5\cot{\theta}}{3\Rey} \right)\partial_x + \frac{5}{6\Rey\Ca}\partial_{xxx} \right]\hat{h} - \left[\frac{5}{2\Rey} + \frac{34}{21}\partial_x \right]\hat{q} + \left[ \frac{1}{3} \right]\hat{f},
  \end{align}
  for the Benney and weighted-residuals models, respectively. We then discretise the problem with finite differences to obtain a system of linear ODEs for each model:
  \begin{equation}
    \label{eqn:disc_benney}\frac{\diff \hat{h}}{\diff t} = (A + BKC)\hat{h},
  \end{equation}
  for Benney and
  \begin{equation}
    \label{eqn:disc_wr}\frac{\diff}{\diff t} \matrix{\hat{h} \\ \hat{q}} = \left[\matrix{A_{hh} & A_{hq} \\ A_{qh} & A_{qq}} + \matrix{B_h \\ B_q} K \matrix{C_h & C_q}\right] \matrix{\hat{h} \\ \hat{q}},
  \end{equation}
  for the weighted residual approach. Here \(A, A_{hh}, A_{hq}, A_{qh}, A_{qq} \in \mathbb{R}^{n\times n}\), \(B, B_h, B_q \in \mathbb{R}^{n\times m}\), \(K \in \mathbb{R}^{m\times p}\), and \(C, C_h, C_q \in \mathbb{R}^{p \times n}\), where there are \(n\) points in the discretisation, \(m\) actuators, and \(p\) observers. We set the discrete cost matrices \(Q = \frac{\beta L}{n}I\) and \(R = I\). We note that there is no particular reason that finite differences has been chosen for the discretisation -- finite elements, finite volumes, and spectral methods would all be equally valid. Having travelled down a hierarchy of increasingly simplified models (summarised in \cref{fig:hierarchy}),
  we are now left with \cref{eqn:disc_benney} and \cref{eqn:disc_wr}, which match the form of \cref{eqn:compact}.

  In the case of perfect interfacial information (\ie where \(C = C_h = I\) and \(C_q = 0\)) previous work~\cite{holroyd2024linear} demonstrated that standard LQR controls, designed for the discretised, linearised, thin-film models, were able to stabilise the flat film from an initially saturated travelling wave. Despite the total lack of predictive capacity at the lowest level of the model hierarchy, given enough actuators (typically at least as many as there were linearly unstable modes), it was possible to control to the target state even up to \(\Rey = 100\).

  \section{Restricted Observations}%
  \label{sec:restricted_observations}
  In computational experiments it is trivial to extract full observations of both the interfacial height \(h\) and the flux \(q\) in both space and time, but in physical experiments this is almost impossible~\cite{heining2013flow}. It is possible to capture observations that are continuous in space (or at the very least high enough resolution that it makes little difference)~\cite{ozgen2002experimental}, and it is often possible to measure at a small number of locations -- for example with lasers~\cite{mouza2000measurement,ishikawa2004measurement} or ultrasound~\cite{wang2019new}.

\begin{figure}[!htb]
    \centering
    \subfloat[Results for the output-feedback LQR method.]{\includegraphics[width=0.48\textwidth]{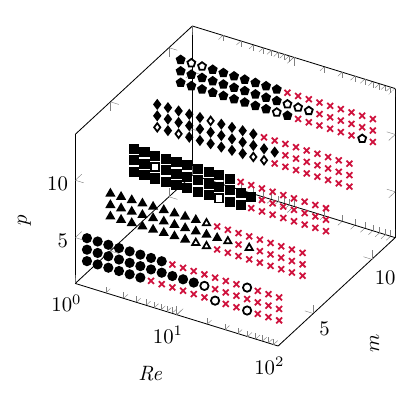}}%
    \hspace{1em}%
    \subfloat[Results for the linear Luenberger-observer method.]{\includegraphics[width=0.48\textwidth]{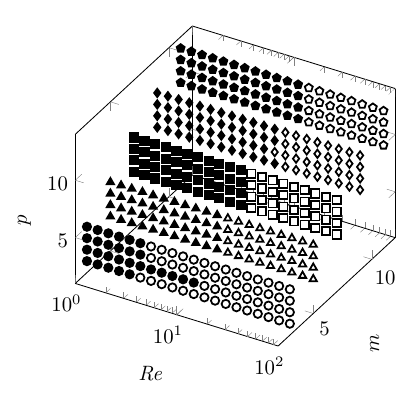}}%
    \caption{In both of the techniques attempted in the authors previous work~\cite{holroyd2024stabilisation} successful control was not possible above \(\Rey = 10\) regardless of the number of actuators \(m\) or observers \(p\). Filled/unfilled shapes correspond to successful/unsuccessful stabilisation. Red crosses denote the failure of the iterative procedure to find a stabilising gain matrix \(K\).}%
    \label{fig:ecc_results}
  \end{figure}

  Unfortunately, even though we have discarded observations of the flux (which is particularly challenging to measure in physical experiments) by setting \(C_q = 0\), the methodology described by the authors~\cite{holroyd2024linear} requires observations of the interface that are continuous in both space and time. In a more recent work they outlined two approaches to handle restricted interfacial observations, where \(p \ll n\). The first of these methods, output-feedback LQR, involves solving the static output feedback problem
  \begin{align}
    \label{eqn:cond1}0 &= A_c\tran P + P A_c + Q + C\tran K\tran R K C, \\
    \label{eqn:cond2}0 &= A_c S + S A_c\tran + I, \\
    \label{eqn:cond3}0 &= RKCSC\tran + B\tran PSC\tran,
  \end{align}
  for the stabilising gain matrix \(K\) (as well as \(P\) and \(S\)), where \(A_c = A+BKC\). When \(C = I\), \cref{eqn:cond1,eqn:cond2,eqn:cond3} collapse to the continuous algebraic Riccati equation (CARE),
  \begin{equation}
    \label{eqn:care}0 = A\tran P + P A_c + - P B R^{-1} B\tran P,
  \end{equation}
  which can be solved directly for \(P`\)~\cite{levine1970determination} or iteratively~\cite{kleinman1968iterative}, and then set \(K = -R^{-1}B\tran P\). The general problem can be solved with expensive iterative methods~\cite{ilka2023novel}, however convergence is only guaranteed from a small set of initial starting matrices, and so, as we can see in \cref{fig:ecc_results}, there is a wide range of the parameter space for which a stabilising gain matrix could not be found. Even when it is possible to compute an optimal \(K\), the control scheme is no longer viable above \(\Rey \approx 10\), representing an order of magnitude loss in performance.

  The second approach was to utilise a Luenberger observer~\cite{luenberger1971introduction}. The idea here is to construct an additional set of ODEs that can be used as a dynamic estimator. Like previous applications of this method to the Kuramoto-Sivashinsky equation~\cite{armaou2000feedback} and the Benney and weighted-residuals models~\cite{thompson2016stabilising}, the first \(p\) linear modes were used to approximate the dynamics, with a forcing term designed to keep the estimator close to the real system. While this solves the numerical challenges inherent to the iterative solver associated with the static output LQR method, as shown in \cref{fig:ecc_results}, it too breaks down at only moderately large Reynolds numbers (\(\Rey \approx 12\)). In the following sections we propose an improved estimator that is better able to capture the behaviour of the Navier-Stokes film.

  \subsection{Observation and Estimation}%
  \label{sub:observation_and_estimation}
  Previous work~\cite{thompson2016stabilising,holroyd2024stabilisation} showcases that, rather than directly passing the observations \(y\) to the gain matrix \(K\), it may be beneficial to instead pass the complete information from the estimator, denoted by \(z\). This means that the actuator magnitudes become \(u = Kz\) rather than \(u = KCx\). Now that \(C\) has been eliminated we can simply compute \(K\) using the standard LQR method as previously developed by the authors~\cite{holroyd2024linear}. However, the question of how to update the estimator remains. The most basic approach is to use the Moore-Penrose inverse to solve \(Cz = y\) in a least-square sense, \ie set \(z = C^{+}y\). This can be viewed as a highly-simplified alternative to the static output feedback LQR problem. Unfortunately since \(p \ll n\) too much information has been lost for this to be successful. The Luenberger observer of the type used in related studies~\cite{thompson2016stabilising,holroyd2024stabilisation} is \(z = \mathcal{F}^{-1}\tilde{z}\), where \(\mathcal{F}^{-1}\) is the inverse Fourier transform (and, if \(p \neq n\), \(\mathcal{F}^{-1}\) also includes interpolation to \(n\) discrete points) and \(\tilde{z} \in \mathbb{C}^p\) is a dynamic estimator of the \(p\) largest modes, which, given a large enough \(p\), should capture the unstable behaviour. However, as we saw in \cref{fig:ecc_results}, this too is eventually insufficient to capture the complex dynamics of the film.

  \begin{figure}[!htb]
    \centering
    \includegraphics[width=\textwidth]{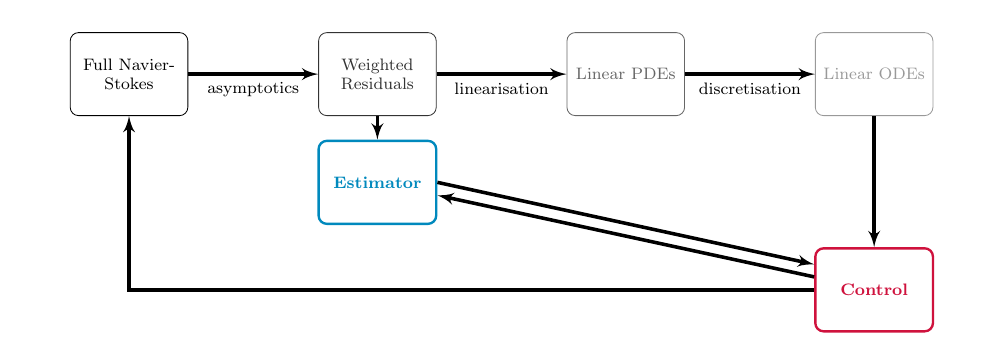}
    \caption{By separating out the estimator model from the control model we can design a better, non-linear estimator to compute the actuator strengths for the full Navier-Stokes simulations, augmenting the approach outlines in \cref{fig:hierarchy}. We note that the full, non-linear Benney equation is not suitable for use as an estimator due to the finite-time blowup.}%
    \label{fig:estimator_hierarchy}
  \end{figure}

  We return to the now augmented hierarchy of model simplifications as depicted in \cref{fig:estimator_hierarchy}. Rather than using a further simplified model of the already linearised, discretised system used to design the controls we now return up the chain to the asymptotic models. Instead of evolving \(z\) using linear dynamics we use the non-linear Benney equation
  \begin{equation}
    \label{eqn:benney_est}z_t = -\left[ \frac{z^3}{3}\left( 2 - 2 z_x \cot\theta + \frac{z_{xxx}}{\Ca} \right) + \Rey \left( \frac{8z^6z_x}{15} - \frac{2z^4f}{3} \right) \right]_x + f + g,
  \end{equation}
  or weighted-residual equations
  \begin{align}
    \label{eqn:wr_est_1}z_t &= f - w_x + g_1, \\
    \label{eqn:wr_est_2}w_t &= -\frac{5w}{2\Rey z^2} + \frac{5z}{6\Rey}\left[ 2 - 2\cot\theta z_x + \frac{z_{xxx}}{\Ca} \right] + \frac{9w^2 z_x}{7z^2} - \frac{17w w_x}{7z} + \frac{f w}{2h} + g_2.
  \end{align}
  Here we distinguish between \(z\), the estimator for the height, and \(w\), the estimator for the flux. The original forcing term remains, but we now add a new forcing \(g\) (or \(g = (g_1, g_2)\) for the coupled weighted-residual model) which is designed to control the system towards the observations. As with finding \(f\) in the main control problem, we are unable to compute \(g\) in a non-linear setting, and so we linearise around the flat state and discretise. Additionally, we set the control term \(g\) to be some linear function of the difference between our observations of the full system \(x\) and mock observations of the estimator \(z\):
  \begin{align}
    \dot{x} &= Ax + BKz, \\
    \dot{z} &= Az + BKz + L(Cx - Cz),
  \end{align}
  where we have set \(f = BKz\). Here we note that \(A\) and \(B\) refer to either the size \(n\) Benney system or the size \(2n\) weighted-residual system -- the linear algebra is identical. Defining the error \(e = x - z\), we obtain a system of ODEs for \(e\):
  \begin{equation}
    \label{eqn:estimator_error}\dot{e} = (A - LC)e.
  \end{equation}
  As with the first LQR problem that we encountered earlier, all that remains is to choose the estimator-gain matrix \(L\) so that \(A - LC\) is stable (\ie all its eigenvalues have negative real part). Noting that eigenvalues of a matrix are invariant under transposition, we can use the LQR algorithm on \((A - LC)\tran = A\tran - C\tran L\tran\) to find the optimal \(L\) (and so we see that the observation problem is the dual of the control problem). This linear forcing term is then substituted back into the non-linear estimator \cref{eqn:benney_est} (or \cref{eqn:wr_est_1,eqn:wr_est_2} for the weighted-residual version) so that we are left with a simplified, but still non-linear estimator with a linear forcing term, as summarised in \cref{fig:estimator_hierarchy}. One concern with this method is that it might impact the computation efficiency of the LQR-based approach, but careful design of the numerical schemes used to update the estimator ensures that this is not the case (see \cref{sec:computational_cost} for additional information).

  Unfortunately, as we saw in \cref{fig:comp_interface}, the singular behaviour of the Benney equation means that it is unsuitable for use as an estimator: the linear forcing term, designed to correct small errors where the behaviour is approximately linear, are unable to arrest the finite-time blowup associated with the travelling wave mode in the Navier-Stokes film. As a result, from this point onwards we will only consider the weighted-residual estimator.

\begin{figure}[!htb]
    \centering
    \includegraphics[width=\textwidth]{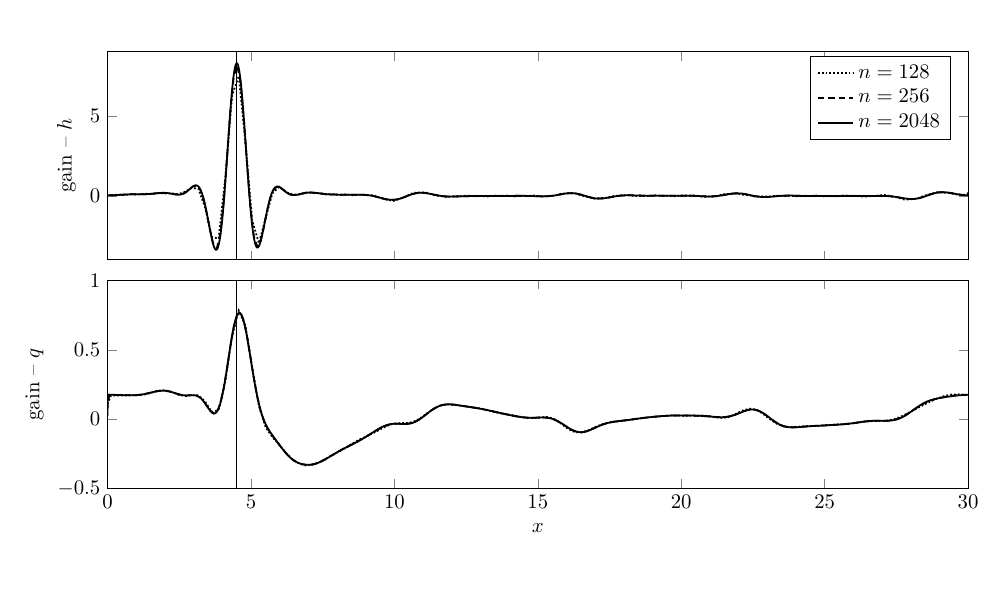}
    \caption{As the number of points in the discretisation -- \(n\) -- is increased, the estimator gain \(L\) converges to the infinite dimensional (but still linear) counterpart, the integral against a continuous weight function. Even at \(n=128\) there is good convergence.}%
    \label{fig:L-vs-N}
  \end{figure}

  Like the original gain matrix \(K\), we should also expect \(L\) to converge to the infinite-dimensional (but still linear) operator that we would recover without the discretisation step~\cite{gibson1983linear}. This is shown in \cref{fig:L-vs-N}.

\begin{figure}[!htb]
    \centering
    \includegraphics[width=\textwidth]{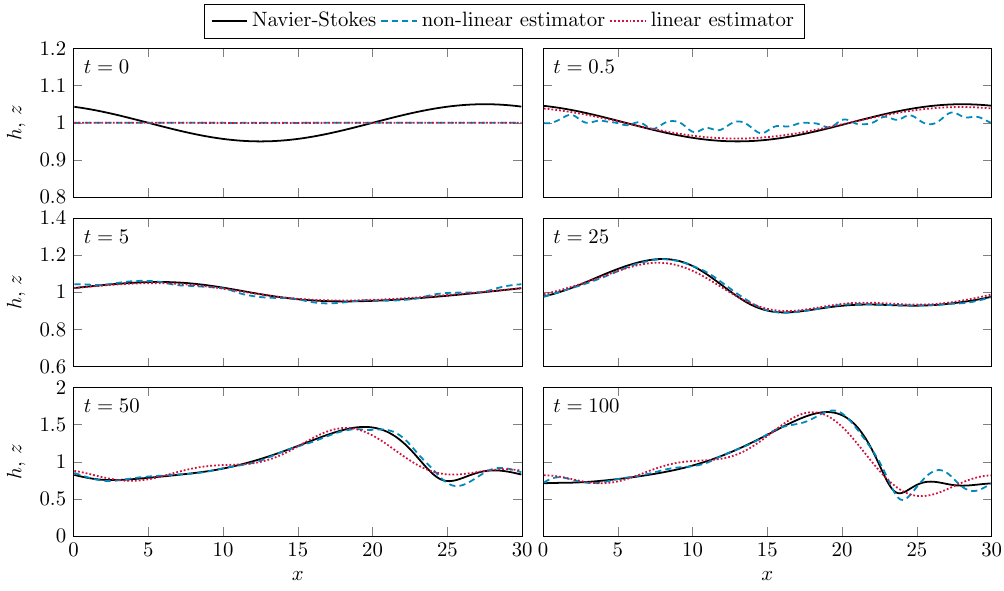}
    \caption{Evolution of the Navier-Stokes film (black), non-linear estimator (blue dashed), and linear estimator (red dotted). Although the linear estimator moves rapidly towards the true profile, once the wave has fully developed the non-linear estimator matches the wave height and position better, as well as capturing the secondary instabilities ahead of the wave.}%
    \label{fig:estimation-snapshots}
  \end{figure}

\begin{figure}[!htb]
    \centering
    \includegraphics[width=\textwidth]{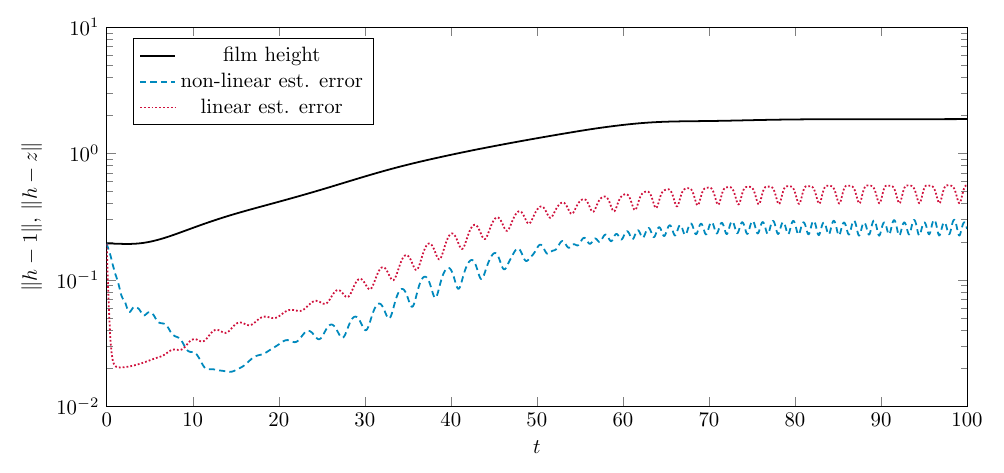}
    \caption{By comparing the 2-norm error between the estimators (linear red dotted, non-linear blue dashed, film perturbation in black) we can see the initial transient behaviour and then the better approximation of the travelling wave by the non-linear estimator.}%
    \label{fig:estimation-performance}
  \end{figure}

  As we can see in \cref{fig:estimation-snapshots,fig:estimation-performance}, the non-linear dynamic estimator is a notable improvement over the linear estimator. The behaviour of the linear estimator is largely driven by the forcing term, resulting in a lag in the peak in \cref{fig:estimation-snapshots}. On the other hand, while the non-linear estimator does suffer from increased oscillations at larger wavenumbers (as seen in the oversized capillary waves in \cref{fig:comp_interface}), the peak height and speed are more accurately captured.
  During the convergence stage we can see similar behaviour in both estimators: an initial move towards the growing interface followed by oscillations as the wave peak passes over the observers. However, once the travelling wave has saturated and the estimator accuracy has plateaued, the frequency of the non-linear estimator's oscillations doubles, and it plateaus at a lower error.

  We also note that, while the eigenvalues for \(A\) and \(A\tran\) (and therefore the number and size of the unstable eigenvalues) are the same, the difficulty of the observation and control problems are different. We start the control problem at the saturated non-linear wave, which is the furthest we can possibly be from the small-perturbation assumption that we made for the linearisation step, and so we should expect the control to be at the limits of its applicability. On the other hand, the estimator starts at \(z(x, 0) = 0\) and evolves as the initial travelling wave develops from a small initial perturbation. This means that the error \(e = h - z\) should always be small, and the corresponding dynamics should therefore be well-approximated by the linear system for which \(L\) is designed. This is one reason why the number of observers required for the estimation problem is lower than the number required for the control problem, as we shall see later.

  \subsection{From Estimation to Control}%
\begin{figure}[!htb]
    \centering
    \includegraphics[width=\textwidth]{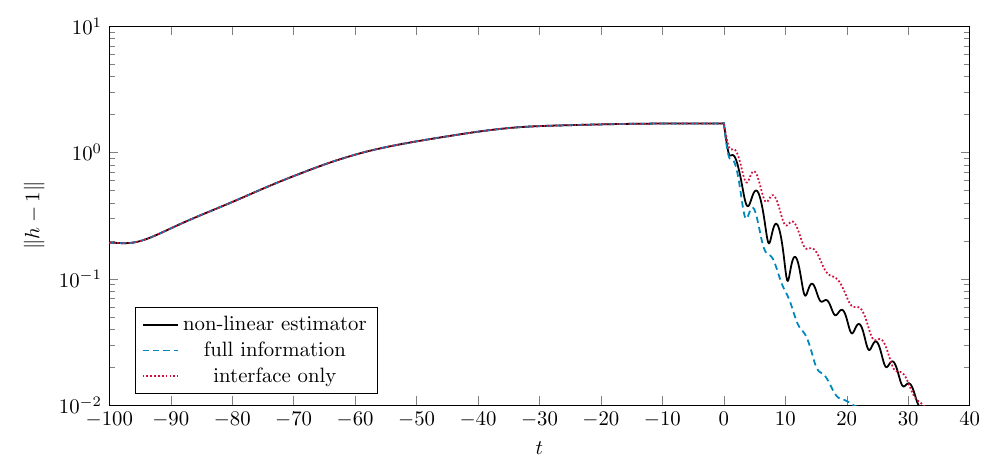}
    \caption{For \(\Rey = 8\) (\(m=5\) and \(p=10\)), LQR controls informed by the non-linear estimator (solid black) perform better than those informed by the full-interface/approximate flux (dotted red), but when given the full flux information (dashed blue) we see additional improvements in the damping rate.}%
    \label{fig:control-performance}
  \end{figure}

  \label{sub:estimation_and_control}
  It is perhaps unsurprising that a non-linear estimator with \(n\) degrees of freedom performs better than a linear estimator with \(p \ll n\) degrees of freedom. In the case of full observations, it is possible to prove that the linearised (but still infinite-dimensional) problem is solvable~\cite{curtain2012introduction} -- \ie that the perturbations from the flat film converge exponentially to zero. The same is true for this method: the linearised problem is provably solvable, and we expect that this translates to some locally stabilisable region in the non-linear case. As \cref{fig:control-performance} shows, the combination of the linear LQR control and non-linear estimator does indeed result in successful exponential damping of the otherwise unstable perturbations.

\begin{figure}[!htb]
    \centering
    \includegraphics[width=\textwidth]{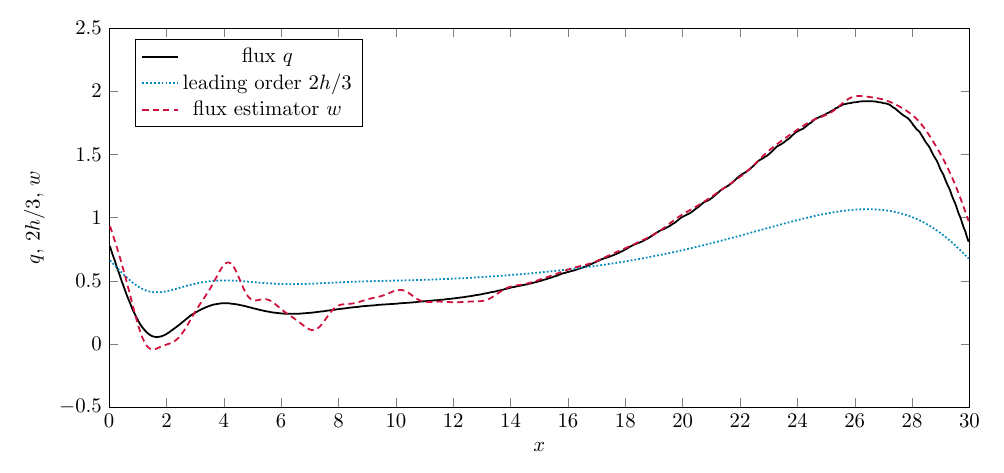}
    \caption{Even for low-to-moderate Reynolds numbers (here \(\Rey = 8\)), the estimator flux \(w\) (dashed red) is much closer to the true flux (solid black) than the leading-order approximation \(q = 2h/3\) (dotted blue). Here \(m=5\) and \(p=10\).}%
    \label{fig:flux-comparision}
  \end{figure}

\begin{figure}[!htb]
    \centering
    \includegraphics[width=\textwidth]{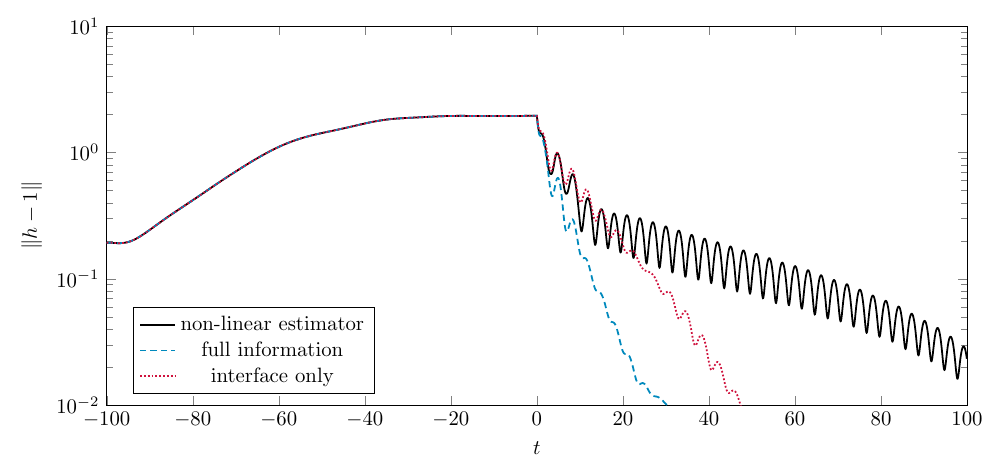}
    \caption{When \(\Rey\) increases to \(12\) (while retaining \(m=5\) and \(p=10\)), the performance of the estimator-informed control deteriorates much faster than the full-information controls.}%
    \label{fig:control-performance-bad}
  \end{figure}

  More interesting is how the controls based on the non-linear estimator compare to the standard LQR controls with full observations. Once a travelling wave has developed (the current authors discussed why this is an appropriate initial condition in~\cite{holroyd2024linear}) and the estimator has stabilised (in \cref{fig:estimation-performance} this happens at around \(t = 75\)), we can activate the controls. \Cref{fig:control-performance} clearly shows that, for low-to-moderate Reynolds numbers, the non-linear estimator achieves the same exponential damping as the original LQR controls. In fact, up to \(\Rey \approx 10\), the estimator is better than the full-interfacial observation LQR control. This is because while the full interfacial information is an advantage, the coupled interface-flux estimator allows for a better approximation of the flux than the leading order method \(q \approx \frac{2}{3}h\) used in \cite{holroyd2024linear,thompson2016stabilising}. This is made obvious in \cref{fig:flux-comparision},
  where the leading-order approximation clearly underestimates the magnitude of the deviation of the flux from the Nusselt flux \(q_N = \frac{2}{3}\) in almost all areas. We can see that when the full-flux information is returned to the LQR control, it then overtakes the estimator in the long run. As \(\Rey\) increases beyond this point and we move further from the asymptotic limit this advantage disappears, as is apparent in \cref{fig:control-performance-bad}.

  As the Reynolds number increases, both the ability of the original LQR control to stabilise the system and the accuracy of the non-linear estimator begin to break down. The performance of the overall control is a combination of the two parts: even though the gain \(K\) might stabilise the Navier-Stokes film for a given Reynolds number, if the information passed to it from the estimator is insufficiently accurate then it will fail, and the gap between \(z\) and \(h\) becomes too large. At this point, the actuators (whose strengths are derived from the now-inaccurate estimator) begin to force different dynamics in \(z\) and \(h\) and the estimator-gain \(LC(h-z)\) is no longer able to overcome the difference. While the linear approximation \cref{eqn:estimator_error} is globally stabilised~\cite{cerpa2014control}, there are no stabilisation guarantees for the non-linear system we are actually working with. For very large Reynolds numbers when the method starts to break down we hypothesise that the introduction of the actuators pushes the estimator system outside of the stabilisable region, which may cause the failure of the controls that we observe.

  With the basic, linear, estimator this failure occurs at \(\Rey = 10\) (\cref{fig:ecc_results}).
  \begin{figure}[!htb]
    \centering
    \includegraphics[width=\textwidth]{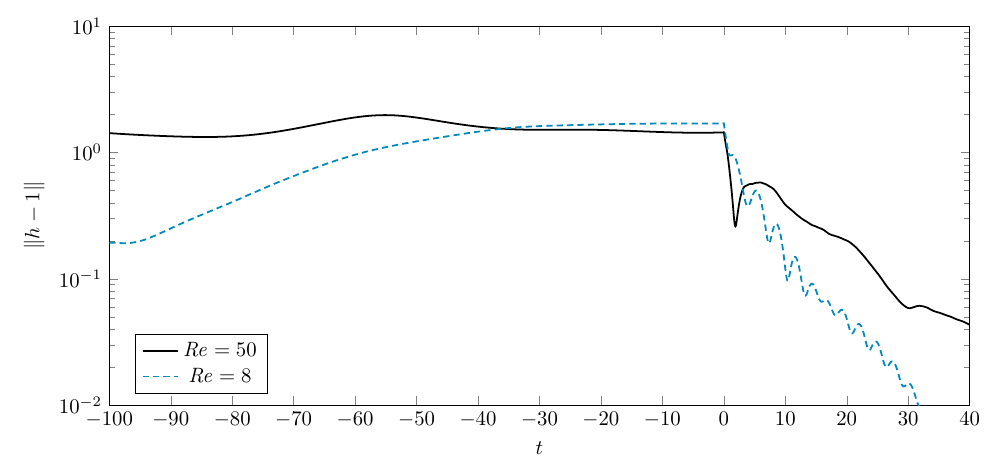}
    \caption{Even at \(\Rey = 50\) (solid black) we are still able to stabilise the interface. As for the full-information case, we must increase the number of actuators (and now observers) to cope with the increased number of unstable modes -- \(m=11\) and \(p=7\). \(\Rey = 8\) case (dashed blue) for comparison.}%
    \label{fig:control-performance-50}
  \end{figure}
  With the non-linear estimator on the other hand, we continue onwards up to \(\Rey = 50\) (\cref{fig:control-performance-50}). It is also worth noting that the number of linearly unstable modes (for both the stabilisation and estimator problems) is
  \begin{equation}
    n_u = 1 + 2 \left\lfloor \frac{L}{2\pi}\sqrt{\Ca \left ( \frac{8}{5}\Rey - 2\cot\theta \right)} \right\rfloor,
  \end{equation}
  which, in this case (\(\Rey = 50\), \(\Ca = 0.01\)), results in \(9\) unstable modes. In previous work~\cite{holroyd2024linear} we showed that, with full information, the film was controllable with \(m > 7\) actuators. Here, despite using only \(p = 7 < n_u = 9\) observers, the estimator follows the behaviour of the true system sufficiently closely such that stabilisation is still achievable.

  \subsection{Observer Placement}%
  \label{sub:observer_placement}
  For the 2D problem with full information, since the boundaries are periodic, it is easy to argue that the \(m\) actuators should be evenly spaced (see section 2.4.2 of \cite{gomes2016control}). Here, when we have \(m = p\) (or \(m \mid p\) or \(p \mid m\)), the same argument suggests that offsetting the actuators and observers by some fixed distance would be optimal. This was discussed by Cimpeanu \etal~\cite{cimpeanu2021active} in the purely proportional control case (\ie when each actuator is set proportionally to the observation from a single observer), and they concluded that a small shift upstream resulted in the fastest rate of damping. A problem more similar to this work -- how to place actuators below a 2D KS equation -- was covered by Tomlin and Gomes~\cite{tomlin2019point}, who concluded that both regular and random arrangements of actuators are able to control the KS equation equally well. Although a full analysis of the interactions of observer and actuator placement are beyond the scope of this work, we have considered the reduced case where we only vary the upstream shift between the observer and actuator locations, \(\delta\), for the case \(m = p = 6\), \(\Rey = 12\).
  \begin{figure}[!htb]
    \centering
    \includegraphics[width=\textwidth]{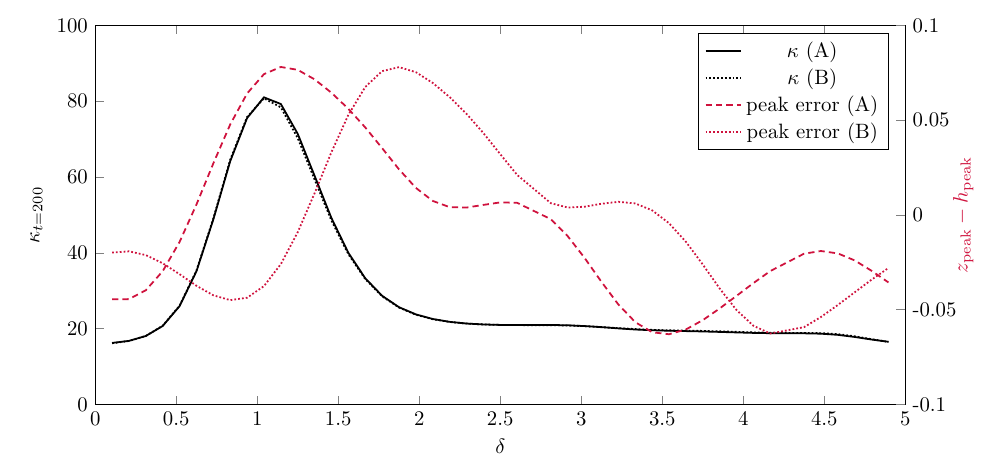}
    \caption{Cost \(\kappa\) at \(t = 200\) as the upstream shift is varied over a single actuator spacing for two different control-activation times (A and B). The cost consistently peaks at \(\delta \approx 1\) regardless of the interfacial position when the controls are activated -- the black lines (solid and dotted) are almost identical. The overestimation of the wave-peak height (red) appears to correlate to the jump in cost for one initial condition (red dashed), but starting at a different wave position (red dotted) shows that the two are not connected.}%
    \label{fig:actuator-placement}
  \end{figure}

  In \cref{fig:actuator-placement} we can see that the final cost (actually the cost at \(t = 200\), when the simulation ended) does indeed vary as we increase the distance that the observers are shifted upstream relative to the actuators, with a strong peak at \(\delta = 1\), representing a small upstream shift of the observers relative to the actuators. Although the jump in cost appeared initially to be correlated to the overestimation of the peak height (the accuracy of the estimator at the peak of the wave oscillates as it passes over the observers -- this is a large component of the oscillations in ~\cref{fig:estimation-performance}), starting the controls at a different time (when the travelling wave position is shifted relative to the observers/actuators and we are at a different position in the oscillations in accuracy) shows that these two are not connected. We leave further investigation into the cause of this phenomenon to a future study.

  \section{Conclusion}%
  \label{sec:conclusion}
  In this work, we extended the basic, linear, Luenberger observer that has been previously shown to be able to stabilise the flat Navier-Stokes film to a non-linear estimator based on a low-dimension asymptotic simplification. We showed that by using this method we can minimise the reduction in the controllable region of the parameter space caused by the dramatic restriction of observations seen in our previous work~\cite{holroyd2024stabilisation}, and that it can even be an improvement over full observations of the interfacial profile for sufficiently low \(\Rey\). We briefly note that the placement of the actuators and observers is no longer a trivial problem, even in the periodic case, and demonstrate that small changes in the relative positions of the actuators and observers can have significant effects.

  There are numerous avenues into which this work could be extended. In addition to the location of actuators and observers, there is the potential for the use of alternative low-dimensional models -- either using the mixed-derivative approach of Ooshida~\cite{ooshida1999surface} or alternative velocity profiles -- both for the control and estimator. There is also the potential for alternative actuation mechanisms such as air-blades (where narrow jets of targeted air are blown at the interface), as well as the coupling of the estimator to non-LQR based approaches such as the electric-field model-predictive controls used by Wray \etal~\cite{wray2022electrostatic}. Current work by the authors includes extending the hierarchical control paradigm to the 3D case, as well as relaxing the assumption of periodic boundaries in favour of a more experimentally realisable inlet-outlet domain.

  \appendix
  \section{Computational Cost}%
  \label{sec:computational_cost}
  One of the principal advantages of LQR-style feedback controls over other types of controls (particularly proportional-integral-derivative controls and model predictive controls) is that the computational complexity is contained in a one-off calculation of the gain matrix at the start of the process -- and since this matrix depends only on the physical parameters it can be reused. The application of the controls depend only on a single matrix-vector multiplication, an \(O(mn)\) operation. The estimator system \cref{eqn:wr_est_1,eqn:wr_est_2} are discretised in space with second-order finite differences (with \(w\) offset from \(z\) by half a grid-spacing to prevent checkerboarding~\cite{sigmund1998numerical}) and stepped forwards in time using a second-order Crank-Nicolson scheme. The implicit component of the Crank-Nicolson step requires solving a system of \(2n\) non-linear equations, which is accomplished using Newton iteration. The na\"ive implementation of this non-linear solve involves the inversion of the associated Jacobian matrix, which is \(2n\)-by-\(2n\), at a cost of \(O(8n^3)\) operations. This quickly becomes overwhelmingly large, far more expensive than solving the full Navier-Stokes equations. Clearly this is impractical -- especially when we consider that in reality these controls need to be applied to physical systems in real time. Fortunately we can utilise the sparse nature of the Jacobian to avoid this prohibitive cost.

  \begin{figure}[!htb]
    \centering
    \includegraphics[width=\textwidth]{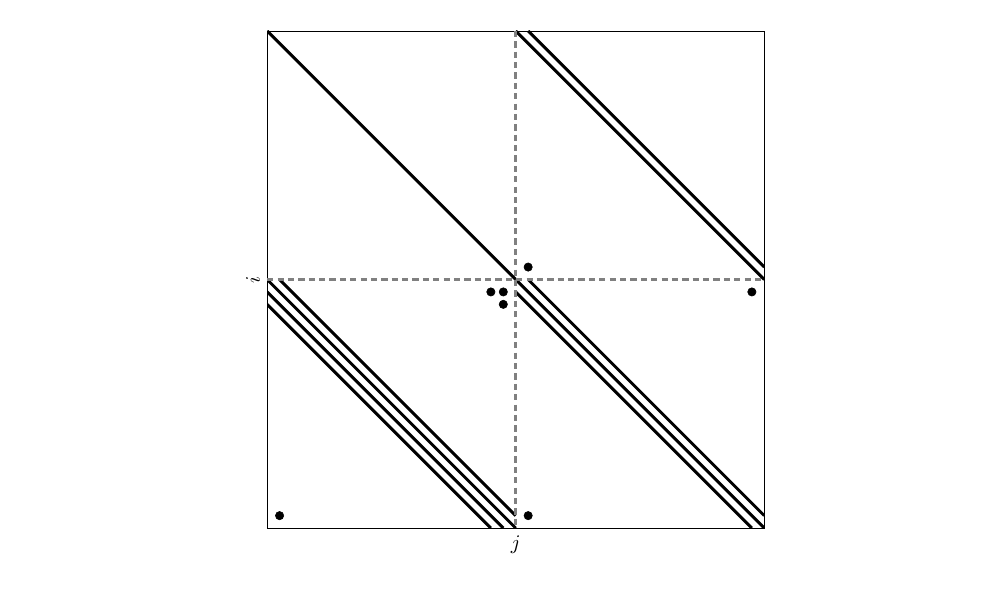}
    \caption{Sparsity pattern for the Jacobian matrix, \cref{eqn:jacobian}. Note how the typical banded structure is augmented by additional individual entries (dots) dues to the periodic boundaries.}%
    \label{fig:jacobian_sparsity}
  \end{figure}

  As we can see in \cref{fig:jacobian_sparsity}, the Jacobian can be divided into four, equally sized, square blocks:
  \begin{equation}
    \label{eqn:jacobian}J = \matrix{J_{zz} & J_{zw} \\ J_{wz} & J_{ww}}
  \end{equation}
  The top left, \(J_{zz}\), is equal to the identity matrix, and so inversion is trivial. Using a blockwise inverse~\cite{lu2002inverses} we can further the problem to the inversion of \(J_{zz}\) and the Schur complement \(S = J_{ww} - J_{wz} J_{zz}^{-1} J_{zw} = J_{ww} - J_{wz} J_{zw}\). \(S\) is a cyclic, circulant, pentadiagonal matrix, which can be constructed and solved in \(O(n)\) steps~\cite{navon1987pent}. The additional cost associated with updating the estimator is therefore only \(O(n)\). Observations indicate that the additional time required in comparison with the standard LQR method is negligible.

  \section{Direct Numerical Simulation}%
  \label{sec:direct_numerical_simulation}
  In lieu of physical experiments, we apply our control schemes to direct numerical simulations of the incompressible Navier-Stokes equations. We perform these simulations using Basilisk~\cite{basilisk2025}, a free, open-source volume-of-fluid (VOF) solver~\cite{scardovelli1999direct} with interface-capturing~\cite{popinet2009accurate,popinet2018numerical}. Basilisk solves for the flow in the liquid and the gas on an adaptive quadtree grid~\cite{popinet2003gerris,popinet2015quadtree}. Although the gas flow has little impact on the liquid, the quadtree grid does not allow for an ``empty'' region. However, the adaptive nature of the grid means that we solve for the gas at a very low resolution, and so it does not have a significant impact on the simulation time. The code used to generate the results in this paper (both for Basilisk DNS and the controls) can be found at \url{https://github.com/OaHolroyd/falling-film-control}.

  \vskip6pt
  \enlargethispage{20pt}

  \ack{Oscar Holroyd is grateful for the computing resources supplied by the University of Warwick Scientific Computing Research Technology Platform (SCRTP) and funding from the UK Engineering and Physical Sciences Research Council (EPSRC) grant EP/S022848/1 for the University of Warwick Centre for Doctoral Training in Modelling of Heterogeneous Systems (HetSys CDT). Radu Cimpeanu and Susana N. Gomes also acknowledge EPSRC Small Grant EP/V051385/1, supporting foundational aspects of this investigation. For the purpose of open access, the authors have applied a Creative Commons Attribution (CC-BY) licence to any arising Author Accepted Manuscript version.}

  \vskip2pc
  \bibliographystyle{RS}
  \bibliography{references}

\end{document}